\documentclass{amsart}
\usepackage{graphicx}
\theoremstyle{definition}

\theoremstyle{remark}

\numberwithin{equation}{section}
\begin{document}

\title{Some gradient estimates for a diffusion equation on Riemannian manifolds}%
\author{Hong Huang}%
\address{School of Mathematics Science,Beijing Normal University,Beijing 100875, P. R. China}%
\email{hhuang@bnu.edu.cn}%

\thanks{Partially supported by NSFC no.10671018.}%

%\date{}%
%\dedicatory{}%
%\commby{}%
% ----------------------------------------------------------------

\begin{abstract}
In this note we present some gradient estimates for the diffusion
equation $\partial_t u=\Delta u-\nabla \phi \cdot \nabla u $ on
Riemannian manifolds, where $\phi $ is a $C^2$ function, which
generalize  estimates of  R. Hamilton's and Qi S. Zhang's on the
heat equation.

\end{abstract} \maketitle
% ----------------------------------------------------------------

\section{Introduction}

Using ideas from Li-Yau [LY], P. Souplet and Qi S. Zhang [SZ] give
an  elliptic estimate for the heat equation on Riemannian manifolds
which is analogous to the Cheng-Yau estimate for the harmonic
functions [CY] and the Hamilton estimates for the heat equation [H].
Later Qi S. Zhang [Z] gives a sharpened local Li-Yau gradient
estimates for the heat equation using ideas from [H] and [LY].

In this note we present some gradient estimates for the diffusion
equation $\partial_t u=\Delta u-\nabla \phi \cdot \nabla u $ on
Riemannian manifolds, where $\phi $ is a $C^2$ function, which
generalize estimates of  R. Hamilton's and Qi S. Zhang's on the heat
equation. Note that similar generalizations of the Cheng-Yau and
Li-Yau estimates have already appeared in X.-D. Li [L], see also Q.
Ruan [R] for a similar generalization ( and for a more general
equation) of the Souplet-Zhang estimate.

\hspace *{0.4cm}

Let $\mathcal M$ be a  Riemannian manifold of dimension $n$, and
$L=\Delta-\nabla \phi \cdot \nabla$, where $\phi $ is a $C^2$
function on $\mathcal M$. Recall the Bakry-Emery Ricci curvature
([BE]) $Ric_{m,n}(L):=Ric+\nabla^2 \phi-\frac {\nabla \phi
\bigotimes \nabla \phi}{m-n}, $  where $m\geq n$ is a constant. (
Here we use the convention that $m=n$ if and only if $L=\Delta$.)
Let $B(x_0,R)$ be a geodesic ball in $\mathcal M$, and
$Q_{R,T}:=B(x_0,R)\times [t_0-T,t_0]$. Then we have the following
generalization of the theorem of  Zhang [Z] mentioned above,

\
\hspace *{0.4cm}

 {\bf Theorem 1.1} Assume that the Bakry-Emery Ricci
curvature $Ric_{m,n}(L) \geq -k$ in $B(x_0,R)$, $k\geq 0$. Suppose
$u$ is any positive solution to the diffusion equation $\partial_t
u=Lu$ in $Q_{R,T}$. Then there exists a dimension constant $c$ such
that

$\frac {|\nabla u|^2}{u^2}-\frac{u_t}{u}\leq
c(\frac{1}{R^2}+\frac{1}{T}+k+(\sqrt{k}+\frac{1}{R})sup \frac
{|\nabla u|}{u})$

in $Q_{R/2,T/2}$.

\hspace *{0.4cm}

The following is a generalization of a theorem due to Hamilton
mentioned above,

\hspace *{0.4cm}

 {\bf Theorem 1.2}
Assume that $\mathcal M$ is compact, and that the Bakry-Emery Ricci
curvature $Ric_{m,n}(L) \geq -k$, $k\geq 0$. Suppose $u$ is any
positive solution to the diffusion equation $\partial_t u=Lu$ in
$\mathcal M \times (0,T]$ ($0<T\leq \infty$) with $u\leq M$. Then

$\frac {|\nabla u|^2}{u^2}\leq (\frac{1}{t}+2k)ln\frac {M}{u}.$

\hspace *{0.4cm}

\section{Proof of  theorems}

 {\bf Proof of theorem 1.1}

As in [Z], we use idea from [H] and [LY].

By  Bakry's generalized Bochner-Weitzenb$\ddot{o}$ck formula ([B]),

$\frac{1}{2}L|\nabla f|^2 = |{\nabla}^2 f|^2+(\nabla Lf)\cdot \nabla
f+(Ric+{\nabla}^2 \phi)(\nabla f,\nabla f),$

 we compute

 $(L-\partial_t)(\frac{|\nabla
u|^2}{u})=\frac{2}{u}|{\nabla}^2 u-\frac{\nabla u \bigotimes \nabla
u}{u}|^2+ 2\frac{(Ric+{\nabla}^2 \phi) (\nabla u, \nabla u)}{u}$

$\geq \frac{2}{nu}(\Delta u -\frac{|\nabla u|^2}{u})^2
+2\frac{(Ric+{\nabla}^2 \phi) (\nabla u, \nabla u)}{u}.$

Note that using the simple inequality $(a+b)^2 \geq
\frac{a^2}{1+\theta}-\frac{b^2}{\theta}$ for any $\theta > 0$, we
have

$(\Delta u -\frac{|\nabla u|^2}{u})^2=(Lu-\frac{|\nabla
u|^2}{u}+\nabla \phi \cdot \nabla u)^2$

$\geq \frac{n}{m}(Lu -\frac{|\nabla u|^2}{u})^2-n\frac {\nabla \phi
\bigotimes \nabla \phi}{m-n}(\nabla u,\nabla u).$

Since $Lu$ is also a solution of the diffusion equation $\partial_t
u=Lu$, by using the assumption $Ric_{m,n}(L) \geq -k$ we have

 $(L-\partial_t)(-Lu+\frac{|\nabla
u|^2}{u})
 \geq\frac{2}{mu}(Lu -\frac{|\nabla
u|^2}{u})^2+\frac{2}{u}(Ric+{\nabla}^2 \phi-\frac {\nabla \phi
\bigotimes \nabla \phi}{m-n})(\nabla u,\nabla u)$

$\geq  \frac{2}{mu}(Lu -\frac{|\nabla u|^2}{u})^2-2k\frac{|\nabla
u|^2}{u}.$

Then, let $H=\frac{|\nabla u|^2}{u^2}-\frac{u_t}{u}$,  we get the
inequality

$(L-\partial_t)H \geq \frac{2}{m}H^2-2k\frac{|\nabla
u|^2}{u^2}-2\nabla H \cdot \nabla ln u.$

Now we use the cut-off function argument as in [LY], [SZ], [Z]. Let
$\psi=\psi(d(x,x_0),t)=\psi(r,t)$ be a smooth cut-off function in
$Q_{R,T}$, such that (i) $0\leq \psi \leq 1$, and $\psi=1$ in
$Q_{R/2,T/4}$; (ii) $\psi$ is decreasing as a radial function in the
spatial variables; (iii)$\frac{|\partial_r \psi|}{\psi^a}\leq
\frac{C_a}{R}$,$\frac{|{\partial_r}^2 \psi|}{\psi^a}\leq
\frac{C_a}{R^2}$ when $0< a < 1$; and (iv)$\frac{|{\partial_t}
\psi|}{\psi^{1/2}}\leq \frac{C}{T}$.

We have

$(L-\partial_t)(\psi H )$

$=\psi(L-\partial_t)H+H(L-\partial_t)\psi+2\nabla \psi \cdot \nabla
H$

 $\geq \frac{2}{m}\psi H^2-2k\psi \frac{|\nabla
u|^2}{u^2}-2\nabla (\psi H)\cdot \nabla ln u+2H\nabla \psi \cdot
\nabla lnu+H(L-\partial_t)\psi+2\frac{\nabla \psi}{\psi}\cdot \nabla
(\psi H)-2\frac{|\nabla \psi|^2}{\psi}H.$

Since $Ric_{m,n}(L)\geq -k$, by the generalized Laplacian comparison
theorem (cf. for example [BQ]) $Lr\leq
(m-1)\sqrt{k}coth(\sqrt{k}r)$, then we have that

$L\psi=(\partial_r \psi) Lr+({\partial_r}^2 \psi )|\nabla r|^2$

$\geq (\partial_r \psi) (m-1)\sqrt{k}coth(\sqrt{k}r)+{\partial_r}^2
\psi$

$\geq (\partial_r \psi) (m-1)(\frac{1}{r}+\sqrt{k})+{\partial_r}^2
\psi.$

Then noticing the properties of $\psi$, we get

$L\psi \geq -c{\psi}^{1/2}(\frac{1}{R^2}+\frac{\sqrt{k}}{R}).$

Suppose the maximum of $\psi H$ is non-negative, otherwise we are
done. Then at the maximum (space-time) point $(y,s)$ of $\psi H$ we
have

$0$

 $\geq \frac{2}{m}\psi H^2-2k\psi \frac{|\nabla
u|^2}{u^2}+2H\nabla \psi \cdot \nabla
lnu-c{\psi}^{1/2}H(\frac{1}{R^2}+\frac{\sqrt{k}}{R}) -H\partial_t
\psi-2\frac{|\nabla \psi|^2}{\psi}H.$

(Note that by Calabi's trick we may assume that $y$ is not in the
cut locus of $x_0$.)

Then using the inequalities (compare with [SZ], [Z])

$-2H\nabla \psi \cdot \nabla lnu$

$\le \frac{\psi H^2}{8m}+c|\frac{\nabla \psi}{\sqrt{\psi}}\cdot
\nabla lnu|^2$

$\le \frac{\psi H^2}{8m}+c(\frac{|\nabla lnu|}{R})^2,$

$H\partial_t \psi$

$\le \frac{\psi H^2}{8m}+c(\frac{\partial_t \psi}{{\psi}^{1/2}})^2$

$\le \frac{\psi H^2}{8m}+c\frac{1}{T^2},$

 and

$2\frac{|\nabla \psi|^2}{\psi}H$

$\le \frac{\psi H^2}{8m}+c(\frac{|\nabla \psi|^2}{{\psi}^{3/2}})^2$

$\le \frac{\psi H^2}{8m}+c\frac{1}{R^4},$

we arrive at

$\psi H^2 \le c (\frac{1}{R^4}+\frac{1}{T^2}+k^2+k(\frac{|\nabla
u|}{u})^2+(\frac{|\nabla lnu|}{R})^2),$

and the result follows.

\hspace *{0.4cm}

{\bf Proof of theorem 1.2}

As in [H], we let $P=\frac{t}{1+2kt}\frac{|\nabla
u|^2}{u}-uln\frac{M}{u}$, then using  Bakry's generalized
Bochner-Weitzenb$\ddot{o}$ck formula, we get

$(L-\partial_t)P\geq 0,$

then the result follows as in [H].

% ----------------------------------------------------------------
\bibliographystyle{amsplain}

\hspace *{0.4cm}

{\bf Reference}

\bibliography{1}[B]Bakry, D.,Un crit$\acute{e}$re de non-explosion pour certaines
diffusions sur une vari$\acute{e}$t$\acute{e}$s riemanniennes
compl$\acute{e}$te, C.R. Acad. Sci. Paris, Ser. I 303(1986), 23-26.

\bibliography{2}[BE]Bakry, D.,Emery, M.,Diffusion hypercontractives,
Sem.Probb. XIX, Lect. Notes in Math. 1123(1985),177-206.

\bibliography{3}[BQ]Bakry, D., Qian, Z.-M., Volume comparison theorems
wi thout Jacobi fields, preprint 2003.

\bibliography{4}[CY] Cheng, S.Y.;Yau, S.T.,Differential equations
on Riemannian manifolds and the geometric applications. Comm. Pure
Appl. Math. 28 (1975), no.3,333-354.

\bibliography{5}[H] Hamilton, R. S., A matrix Harnack estimate for
the heat equation. Comm. Anal. Geom. 1(1993), no.1, 113-126.

\bibliography{6}[LY] Li, P.;Yau, S.T., On the parabolic kernel of
the Schr$\ddot{o}$dinger operator, Acta Math. 156(1986), 153-201.

\bibliography{7}[L] Li,X.-D., Liouville theorems for symmetric
diffusion operators on complete Riemannian manifolds, J. Math. Pure
Appl. 84(2005),1295-1361.

\bibliography{8}[R] Ruan, Q., An elliptic type gradient estimate for
the Schr$\ddot{o}$dinger equation, arXiv:0705.4178.

\bibliography{9}[SZ]Souplet, S., Zhang, Qi S., Sharp gradient
estimates and Yau's Liouville theorem for the heat equation on
noncompact manifolds, Bull. Lond. Math. Soc.38(2006),no.6,
1045-1053.

\bibliography{10}[Z]Zhang, Qi S.,Some gradient estimates for the heat
equation on domains and for an equation by Perelman, Inter. Math.
Res. Notices 2006.

\end{document}